\documentstyle{article}

\newcommand{\U}{{\bf U}}

\newcommand{\QED}{\hfill$\Box$\medskip}

\begin{document}

\title{About Stability of Irreducibility for Germs of Holomorphic Functions}
\author{\and Huayi Zeng\thanks{Thanks for helpful discussion with Sorin Popescu, Yusuf Mustopa, and Luis E. Lopez} \\ State University of New
York at Stony Brook \\ (hzeng@math.sunysb.edu)}

\date{April 11, 2005}
\maketitle
\begin{center}
{\bf Abstract}
\end{center}

This survey is about irreducibility for germs of a holomorphic
functions $f$. I will show that when the dimension of the domain
$U$ of this holomorphic function $f$ is greater than 2, the
irreducibility of germs are not necessary to be stable. That
means, if the germ of $f$ at point $p$ is irreducible in the stalk
of holomorphic functions at $p$, this does NOT means there exists
an open neighborhood $V\subset U$ of this point $p$, such that for
any point $q\in V$, the germ of $f$ at $q$ is irreducible at the
stalk of holomorphic functions at $q$

\section{Introduction}

Let $U$ be an open set in \textbf{C}$^n$ which contains 0, $f$ be
a holomorphic function defined on $U$, $f_p$ is the germ of $f$ at
point $p\in U$.\\

For any two holomorphic functions $g,h$ defined on $U$,if
$g_0,h_0$ are relatively prime with each other, then with the help
of resultants, we know that $g,h$ are relatively prime with each
other nearby. Precisely to say, that means their exists an open
neighborhood $V\subset U$ of 0, such that for any point $q\in V$,
$g_q$ and $h_q$ are relatively prime with each other. In this
sense, we can say that \textbf{Being co-prime is a stable
property.}\\

Can we say \textbf{Irreducibility is a stable property}?In the
case of dimension 2, the answer is positive, and the proof is
easy. But in the case of dimension 3, I will present a polynomial
as counter-example.

\section{Proof for the Case of Dimension 2}
\textbf{Statement}: For any holomorphic function $f=f(z_1,z_2)$ on
$U\subset$ \textbf{C}$^2$($0\in U$), and the germ of f at origin
is irreducible, then their exists an open neighborhood $V\subset
U$ of 0, such that for any point $q\in V$, $f_q$ is
irreducible.(\textbf{Remark}:If $f(p)\neq 0$, the $f$ is irreducible at $p$. So we only need to care about zero points of $f$.)\\\\
\textbf{Proof:} Without the loss of generality, we can assume
$f(0,z_2)$ is not identically 0 near the origin, and $f(0,0)=0$.\\\\
let $w=z_2^d+e_1(z_1)z_2^{d-1}+\cdots+e_{d-1}(z_1)+e_d(z_1)$ be a
Weierstrass polynomial of $f$ near 0.\\\\
Because $w$ is irreducible at 0, so $w$ and $\frac{\partial
w}{\partial z_2}$ are relatively prime near 0. Then the resultant
of $w$ and $\frac{\partial w}{\partial z_2}$ is not zero. Then the
common zero loci of $w$ and $\frac{\partial w}{\partial z_2}$ are
discrete near 0.\\\\
From above, we know that their exists an open set $V(0\in V\subset
U)$, such that in $\U$, (0,0) is the only zero point of $w$ which
is POSSIBLE to be singular.(since for other points in $q\in U$,
$\frac{\partial w}{\partial z_2}(p)\neq 0$ ).We can conclude that
at any zero point  $p(p\neq 0)$ of $w$ in $V$,$w$ is a local
complex parameter near $p$. Since $w$ is a local complex parameter
near $p$, then the germ of $w$ at $p$ is irreducible.\\\\
Finally, because $w$ is a Weierstrass polynomial of $f$ at 0, then
we know that in $V$, the irreducibility of $f$ is as the same as t
that of $w$.\QED

\section{A Counter Example in Dimension 3} \label{se:C0}
In the case of dimension 3, the statement should be:\\\\
\textbf{Statement}: For any holomorphic function
$f=f(z_1,z_2,z_3)$ on $U\subset$ \textbf{C}$^3$($0\in U$), and the
germ of f at origin is irreducible, then their exists an open
neighborhood $V\subset U$ of 0, such that for any point $q\in V$,
$f_q$ is irreducible.\\\\
But unfortunately, this statement is not true.In this section, I
will present, a polynomial of three variables,
as a counter example.\\\\
This polynomial is $f=z_3^2-z_1z_2^2$.\\\\

\subsection{Irreducibility of $f$ at origin}
Obviously, near 0, $f$ is a Weierstrass polynomial of itself(we
choose $z_3$ as the polynomial variable).Now, we will show the
irreducibility at origin by means of contradiction. \\\\
If $f$ is not irreducible at origin, then its Weierstrass
polynomial is decomposable at origin as a Weierstrass
Polynomial.Assume that,near origin,
$f=(z_3-g_(z_1,z_2))(z_3-h(z_1,z_2))$, here $g,h$ are holomorphic
functions of variable $z_1,z_2$ near 0, and g(0,0)=h(0,0)=0.\\\\
From the factorization $f=(z_3-g_(z_1,z_2))(z_3-h(z_1,z_2))$, we
know that $g+h=0,gh=-z_1z_2^2$, which implies $g^2=z_1z_2^2$ near
0.\\\\
But if $g^2=z_1z_2^2$ near 0. Then for some $\varepsilon\in$
\textbf{C} whose norm is small enough,
$g^2(z_1,\varepsilon)=\varepsilon ^2 z_1$ near 0. But just from
elementary knowledge of functions of one complex variable, we know
this is not possible.\\\\ From argument above, we know $f$ is
irreducible at origin.

\subsection{Further Argument}
At point $p=(z,0,0)(z\neq 0)$, we know that $f(p)=0$, and easily
we can factorize $f$ as $f=(z_3+z_2r)(z_3-z_2r)$ near $p$, here
$r$ is a one-variable holomorphic function such that $r^2=z_1$
near $(z,0,0)$(Because z is not 0, so we can take square-root of
$z_1$ near by.).\\\\
From the argument in \textbf{3.2}, we know that, in any
neighborhood $U$ of origin, there EXISTS some point $p$ such that
$f$ is not irreducible at $p$. This fact can destroy our statement
at the beginning of this section.

\end{document}